\documentclass[12pt]{article}
\usepackage{amsmath,latexsym,amssymb}
\usepackage{amsthm}
\usepackage{epsfig,graphicx}

\newcommand\NN{\mathbb{N}}
\newcommand\RR{\mathbb{R}}
\newcommand\ZZ{\mathbb{Z}}

\newcommand\BB{\mathcal{B}}
\newcommand\CC{\mathcal{C}}
\newcommand\DD{\mathcal{D}}

\newcommand\GG{\mathcal{G}}
\newcommand\HH{\mathcal{H}}

\newcommand\eps{{\varepsilon}}
\DeclareMathOperator\sP{P}   
\newcommand{\rP}{\mathrm{P}} 
\DeclareMathOperator\sE{E}   
\newcommand{\rE}{\mathrm{E}} 

\DeclareMathOperator\dist{dist}

\DeclareMathOperator\diam{diam}
\DeclareMathOperator\defin{def}
\DeclareMathOperator\height{height}
\renewcommand{\mod}{\operatorname{mod}}

\newtheorem{theorem}{Theorem}[section]
\newtheorem{definition}[theorem]{Definition}
\newtheorem{corollary}[theorem]{Corollary}

\newtheorem{lemma}[theorem]{Lemma}
\newtheorem{remark}[theorem]{Remark}
\newtheorem{proposition}[theorem]{Proposition}

\newcommand{\address}{Address: Department of Mathematics, University of North Texas, 1155 Union Circle \#311430, Denton, TX 76203-5017, USA; E-mail: allaart@unt.edu}

\title{How large are the level sets of the Takagi function?}
\author{Pieter C. Allaart \footnote{\address}}

\begin{document}

\maketitle

\begin{abstract}

Let $T$ be Takagi's continuous but nowhere-differentiable function. 
This paper considers the size of the level sets of $T$ both from a probabilistic point of view and from the perspective of Baire category. We first give more elementary proofs of three recently published results. The first, due to Z. Buczolich, states that almost all level sets (with respect to Lebesgue measure on the range of $T$) are finite. The second, due to J. Lagarias and Z. Maddock, states that the average number of points in a level set is infinite. The third result, also due to Lagarias and Maddock, states that the average number of local level sets contained in a level set is $3/2$. In the second part of the paper it is shown that, in contrast to the above results, the set of ordinates $y$ with uncountably infinite level sets is residual, and a fairly explicit description of this set is given. In addition, it is shown that most level sets (in the sense of Baire category) contain infinitely many local level sets, and that a continuum of level sets even contain uncountably many local level sets. Finally, most of the main results are extended to a somewhat more general family of nowhere-differentiable functions. 


\end{abstract}

\section{Introduction}

Takagi's continuous nowhere differentiable function is defined by
\begin{equation}
T(x)=\sum_{n=0}^\infty \frac{1}{2^n}\phi(2^n x),
\label{eq:Takagi-def}
\end{equation}
where $\phi(x)=\dist(x,\ZZ)$, the distance from $x$ to the nearest integer. Since its initial discovery in 1903 by Takagi \cite{Takagi} and subsequent rediscovery by Van der Waerden \cite{vdW}, Hildebrandt \cite{Hildebrandt} and others, much has been written about this function. Arguably the simplest proof of the nowhere-differentiability of $T$ was given by Billingsley \cite{Billingsley}, and his argument was adapted by Cater \cite{Cater} to show that $T$ does not possess a finite {\em one-sided} derivative at any point. In 1959, Kahane \cite{Kahane} showed that the maximum value of $T$ is $\frac23$, and the set of points where this is attained is a Cantor set of Hausdorff dimension $\frac12$. Later still, with the increased popularity of fractals following B. Mandelbrot's work, $T$ became known as an example of a self-affine function, although its Hausdorff dimension of $1$ classifies it as a ``borderline fractal". Many other properties of the Takagi function and various of its generalizations have been investigated. Nonetheless, the Takagi function itself has been slow to give up some of its deepest secrets. For example, it has only recently been established at which set of points $T(x)$ has an infinite derivative (Allaart and Kawamura \cite{AK}, Kr\"uppel \cite{Kruppel2}), and there are still many open questions regarding the level sets of $T$. This paper aims to answer some of these questions. It complements recent work on the level sets of the Takagi function by Knuth \cite{Knuth}, Buczolich \cite{Buczolich}, Maddock \cite{Maddock}, and Lagarias and Maddock \cite{LagMad1,LagMad2}.

For $y\in[0,\frac23]$, define 
$$L(y)=\{x\in[0,1]: T(x)=y\}.$$
Thus, $L(y)$ is the level set at level $y$ of the Takagi function. Since $T(x)>0$ for all $x\in(0,1)$, the simplest level set is $L(0)=\{0,1\}$. At the other extreme, Kahane \cite{Kahane} showed that $L(\frac23)$ is the set of all $x\in[0,1]$ whose binary expansion $x=0.b_1 b_2 b_3\dots$ satisfies $b_{2i-1}+b_i=1$ for all $i\in\NN$. This is equivalent to saying that the quarternary expansion of $x$ contains only $1$'s and $2$'s. As a result, $L(\frac23)$ is a Cantor set of Hausdorff dimension $\frac12$. Surprisingly, a more general study of the level sets of $T$ was apparently not undertaken until 2005, when Knuth \cite[p.~103]{Knuth} published an algorithm for determining $L(y)$ for rational $y$. A few years later, Buczolich \cite{Buczolich} showed that almost all level sets (with respect to Lebesgue measure) are finite. Shortly afterwards, Maddock \cite{Maddock} proved that the Hausdorff dimension of any level set of $T$ is at most $0.668$, and conjectured an upper bound of $\frac12$; his conjecture was recently proved by de Amo et al.~\cite{ABDF}. Lagarias and Maddock \cite{LagMad1,LagMad2} introduced the concept of a {\em local level set} to prove a number of new results. For instance, they show that the {\em average} cardinality of all level sets is infinite, and that the set of ordinates $y$ for which $L(y)$ has strictly positive Hausdorff dimension is of full Hausdorff dimension 1. Combined with the result of Buczolich, these results sketch a complex picture of the totality of level sets of the Takagi function. 

The aim of the present paper is to give more direct proofs of some of the above-mentioned results, and to explore the sizes of the level sets in more detail, both from the probabilistic perspective and from the point of view of Baire category. The results can be summarized by saying that most level sets are finite when viewed probabilistically, but most are uncountably infinite when viewed through the lens of Baire category. 

First, in Section \ref{sec:prelim}, we state the well-known fact that the graph of $T$ contains everywhere miniature copies of itself, which we call {\em humps}. These humps are classified both by their generation and their order (or size), and for many results in the paper it is important to count carefully the number of humps of a given order, either overall or in a given generation. This counting involves the famous Catalan numbers.

Section \ref{sec:probability} views the level sets from a probabilistic perspective. The results in this section are not new, but we show that they can all be derived with little effort from a single key observation: Lemma \ref{lem:bijection} below shows that when all the humps of generation $1$ are removed from the graph of $T$, what remains is a set which intersects every horizontal line between $y=0$ and $y=\frac12$ in exactly two points. The ideas of the proof are implicit in the proof of Buczolich \cite[Theorem 9]{Buczolich}, but we identify them more explicitly here to show how much more can be done with them.
We use the lemma first to derive Buczolich's result that almost all level sets are finite (Theorem \ref{thm:finite-ae}). Then, in Theorem \ref{thm:infinite-expectation}, we use it to prove a result of Lagarias and Maddock \cite{LagMad2}, which states that the average cardinality of the level sets is infinite. 
The last theorem in this section involves the notion of {\em local level set}, introduced by Lagarias and Maddock \cite{LagMad1}. A local level set is a set of abscissas $x$ which can be obtained from one another by certain combinatorial operations on their binary expansions that leave the value of $T(x)$ unchanged -- see Section \ref{subsec:local} for a precise definition. We give an alternative and more elementary proof of one of their results, namely that the average number of local level sets contained in a randomly chosen level set is $\frac32$ (Theorem \ref{thm:local-level-sets}). Our proof takes full advantage of the fact that almost all level sets are finite.

Section \ref{sec:Baire} examines the level sets from the point of view of Baire category, and the results there contrast sharply with those of Section \ref{sec:probability}. To the best of the author's knowledge, the results of this section are new. We first show in Theorem \ref{thm:topology} that, for each finite number $m$, the set of ordinates $y$ for which $|L(y)|\leq m$ is nowhere-dense, and hence the set $\{y\in[0,\frac23]: |L(y)|<\infty\}$ is of the first category. This basic result is then refined in several ways. For example, in Theorem \ref{thm:S-uncountable} the set $S_\infty^{uc}$ of ordinates $y$ with {\em uncountably} infinite level sets is shown to be residual, and a fairly explicit description of this set is given. This includes the observation that $S_\infty^{uc}$ does not contain any dyadic rational points. To end the paper, we show in Theorem \ref{thm:infinite-local-level-sets} that most level sets (in the sense of Baire category) contain infinitely many {\em local} level sets, and that a continuum of level sets even contain uncountably many local level sets. Seen in this light, the last result of Section \ref{sec:probability} appears quite remarkable.

In Section \ref{sec:general} we generalize several of the main results to a version of the Takagi function in which the summands are multiplied by arbitrary signs. The bottom line is the same: almost all level sets are finite, but the average cardinality of the level sets is infinite, and the typical level set is uncountable. The idea of a local level set remains relevant, and we show that the average number of local level sets contained in a level set is between $\frac32$ and $2$.

In a separate paper \cite{Allaart2}, we focus on the cardinalities of the finite level sets of $T$. One of the main results of that paper is that any even positive integer is the cardinality of some level set. We show that more than 60\% of all level sets (in the sense of Lebesgue measure) have exactly two elements, and examine which other cardinalities occur with positive probability when an ordinate $y$ is chosen at random. However, this question is still not fully resolved.

In recent years, the Takagi function has appeared in a surprising range of applications. For instance, Trollope \cite{Trollope} and Delange \cite{Delange} discovered that $T$ was the missing piece of the puzzle in the binary digital sum problem in number theory. The Takagi function also arises naturally in the limit in certain counting problems in graph theory; see Frankl et al. \cite{Frankl} or Guu \cite{Guu}. It plays a role in analyzing the attractors of certain chaotic dynamical systems (see Yamaguti et al. \cite{Yamaguti}), and can even be used to state an equivalent formulation of the Riemann hypothesis (see Balasubramanian et al. \cite{Balasub}). More recently, Tabor and Tabor \cite{Tabor} reported that $T$ is the extremal case in a study of approximate convexity. This wide range of different contexts in which the Takagi function appears makes it likely that almost any intrinsic aspect of this function will find some use, thereby justifying its continued study.

\section{Preliminaries} \label{sec:prelim}

In this paper, $|.|$ will always denote cardinality; the diameter of a set $A$ will be denoted by $\diam(A)$. We denote by $\ZZ_+$ the set of nonnegative integers, by $\NN$ the set of positive integers, and by $\pi_Y(A)$ the projection of a set $A\subset \RR^2$ onto the $y$-axis.

We first recall some known facts about the Takagi function, and introduce important notation and terminology. One of the most important aspects of $T$ for the purposes of this paper is its symmetry with respect to $x=\frac12$:
\begin{equation*}
T(1-x)=T(x) \qquad\mbox{for all $x\in[0,1]$}.
\end{equation*}

Next, define the {\em partial Takagi functions}
\begin{equation*}
T_k(x):=\sum_{n=0}^{k-1}\frac{1}{2^n}\phi(2^n x), \qquad k=1,2,\dots.
\end{equation*}
Each function $T_k$ is piecewise linear with integer slopes. In fact, the slope of $T_k$ at a non-dyadic point $x$ is easily expressed in terms of the binary expansion of $x$. We define the binary expansion of $x\in[0,1)$ by
\begin{equation}
x=\sum_{n=1}^\infty \frac{\eps_n}{2^n}=0.\eps_1\eps_2\dots\eps_n\dots, \qquad\eps_n\in\{0,1\},
\label{eq:binary-expansion}
\end{equation}
with the convention that if $x$ is dyadic rational, we choose the representation ending in all zeros. For $k=0,1,2,\dots$, let
\begin{equation*}
D_k(x):=\sum_{j=1}^k(1-2\eps_j)=\sum_{j=1}^k(-1)^{\eps_j}
\end{equation*}
denote the excess of $0$ digits over $1$ digits in the first $k$ binary digits of $x$. Then it follows directly from  
\eqref{eq:Takagi-def} that the slope of $T_k$ at a non-dyadic point $x$ is $D_k(x)$. The Takagi function itself can be expressed in terms of the sequence $\{D_n(x)\}$ via the formula
\begin{equation}
T(x)=\frac12-\frac14\sum_{n=1}^\infty(-1)^{\eps_{n+1}(x)}\frac{D_n(x)}{2^n},
\label{eq:D-expression}
\end{equation}
as shown by Lagarias and Maddock \cite[Section 2]{LagMad1}; see also Lemma \ref{lem:general-D-expression} below. This formula yields an easy proof of the following important fact.

\begin{lemma} \label{lem:equivalence}
If $|D_n(x)|=|D_n(x')|$ for every $n$, then $T(x)=T(x')$.
\end{lemma}

\begin{proof}
Let $Z=\{n\geq 0: D_n(x)=0\}$, and enumerate the elements of $Z$ as $0=n_0<n_1<n_2<\dots$. If $D_n(x)=-D_n(x')$ for $n_i<n<n_{i+1}$ (where possibly $n_{i+1}=\infty$), then it must be the case that $\eps_n(x)=1-\eps_n(x')$ for $n_i<n\leq n_{i+1}$. Since $D_{n_{i+1}}(x)=D_{n_{i+1}}(x')=0$, it follows that the part of the summation in \eqref{eq:D-expression} over $n_i<n\leq n_{i+1}$ is the same for $x$ as for $x'$. This can be repeated for all intervals $(n_i,n_{i+1}]$ on which $D_n(x)\neq D_n(x')$. Hence, $T(x)=T(x')$.
\end{proof}

Later in this section we will sketch an alternative proof which does not rely on the formula \eqref{eq:D-expression}, but which exploits more intuitively the symmetry and self-similarity properties of the graph of $T$.
An important consequence of Lemma \ref{lem:equivalence} is that, if $D_n(x)=0$ for infinitely many indices $n$, then the level set $L(T(x))$ contains a Cantor set, as there are $2^{\aleph_0}$ many points $x'$ such that $|D_n(x)|=|D_n(x')|$ for each $n$.

The next lemma is well known (e.g. Kahane \cite{Kahane}), but a short proof is included in order to keep this paper self-contained.

The term `balanced' in the following definition is taken from Lagarias and Maddock \cite{LagMad1}.

\begin{definition} \label{def:balanced}
{\rm
A dyadic rational of the form $x=0.\eps_1\eps_2\dots\eps_{2m}$ is called {\em balanced} if $D_{2m}(x)=0$. If there are exactly $n$ indices $1\leq j\leq 2m$ such that $D_j(x)=0$, we say $x$ is a balanced dyadic rational of {\em generation} $n$. By convention, we consider $x=0$ to be a balanced dyadic rational of generation $0$.

The set of all balanced dyadic rationals is denoted by $\BB$. For each $n\in\ZZ_+$, the set of balanced dyadic rationals of generation $n$ is denoted by $\BB_n$. Thus, $\BB=\bigcup_{n=0}^\infty \BB_n$.
}
\end{definition}

The following lemma states in a precise way that the graph of $T$ contains everywhere small-scale similar copies of itself. Let
\begin{equation*}
\mathcal{G}_T:=\{(x,T(x)): 0\leq x\leq 1\}
\end{equation*}
denote the graph of $T$ over the unit interval $[0,1]$.

\begin{lemma} \label{lem:similar-copies}
Let $m\in\NN$, and let $x_0=k/2^{2m}=0.\eps_1\eps_2\dots\eps_{2m}$ be a balanced dyadic rational. Then for $x\in[k/2^{2m},(k+1)/2^{2m}]$, we have
\begin{equation*}
T(x)=T(x_0)+\frac{1}{2^{2m}}T\left(2^{2m}(x-x_0)\right).
\end{equation*}
In other words, the part of the graph of $T$ above the interval $[k/2^{2m},(k+1)/2^{2m}]$ is a similar copy of the full graph $\mathcal{G}_T$, reduced by a factor $1/2^{2m}$ and shifted up by $T(x_0)$.
\end{lemma}

\begin{proof}
This follows immediately from the definition \eqref{eq:Takagi-def}, since the slope of $T_{2m}$ over the interval $[k/2^{2m},(k+1)/2^{2m}]$ is equal to $D_{2m}(x_0)=0$, and $T(x_0)=T_{2m}(x_0)$.
\end{proof}

\begin{figure}
\begin{center}
\epsfig{file=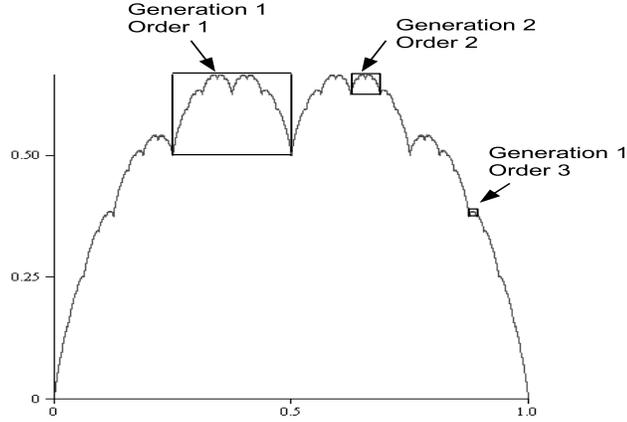, height=.35\textheight, width=.7\textwidth}
\caption{The graph of $T$, with humps of various orders and generations highlighted. The rectangles shown are, from left to right, $K(1/4)$, $K(5/8)$ and $K(7/8)$. Note that in binary, $1/4=0.01$, $5/8=0.1010$, and $7/8=0.111000$.}
\label{fig:humps}
\end{center}
\end{figure}

\begin{definition} \label{def:humps}
{\rm
For a balanced dyadic rational $x_0=k/2^{2m}$ as in Lemma \ref{lem:similar-copies}, define 
\begin{gather*}
I(x_0)=[k/2^{2m},(k+1)/2^{2m}], \qquad J(x_0)=T(I(x_0)),\\
K(x_0)=I(x_0)\times J(x_0),\\
H(x_0)=\mathcal{G}_T\cap K(x_0).
\end{gather*}
By Lemma \ref{lem:similar-copies}, $H(x_0)$ is a similar copy of the full graph $\mathcal{G}_T$; we call it a {\em hump}. Its height is $\diam(J(x_0))=\frac23{(\frac14)}^m$, and we call $m$ its {\em order}. By the {\em generation} of the hump $H(x_0)$ we mean the generation of the balanced dyadic rational $x_0$. A hump of generation $1$ will be called a {\em first-generation hump}. By convention, the graph $\mathcal{G}_T$ itself is a hump of generation $0$. If $D_j(x_0)\geq 0$ for every $j\leq 2m$, we call $H(x_0)$ a {\em leading hump}. See Figure \ref{fig:humps} for an illustration of these concepts.

We denote by $\HH=\{H(x_0): x_0\in\BB\}$ the set of all humps, and by $\HH'$ the subset of $\HH$ consisting of all leading humps.
}
\end{definition}

It is clear that each hump of generation $n$ is a subset of a hump of generation $n-1$. Note also that the projections of the first-generation humps onto the $x$-axis have disjoint interiors. 

We can now sketch a more geometric proof of Lemma \ref{lem:equivalence}. 

\begin{proof}[Second proof of Lemma \ref{lem:equivalence} (sketch)]
Let $x\in[0,1]$, and assume that $x$ is not dyadic rational. Given a fixed number $m$ such that $D_{2m}(x)=0$, there is a unique point $x'$ such that $D_n(x')=D_n(x)$ for all $n\leq 2m$, but $D_n(x')=-D_n(x)$ for all $n>2m$. Note that $x\in I(x_0)$ with $x_0=k/2^{2m}$ for some integer $k$. It is easy to see that $x-k/2^{2m}=(k+1)/2^{2m}-x'$. By Lemma \ref{lem:similar-copies}, the graph of $T$ above $I(x_0)$ is a hump, and is hence symmetric about its central vertical axis. Therefore, $T(x)=T(x')$. If $x$ and $x'$ are any two non-dyadic points satisfying the hypothesis of Lemma \ref{lem:equivalence}, $x'$ can be obtained from $x$ by countably many operations of the above type, so a simple application of the continuity of $T$ gives that $T(x)=T(x')$. Clearly, if the lemma holds for non-dyadic points, it holds for dyadic points as well, again by continuity of $T$.
\end{proof}

Many of the results in this paper depend on a careful count of humps of a given order. This involves the 
{\em Catalan numbers}
\begin{equation*}
C_n:=\frac{1}{n+1}\binom{2n}{n}, \qquad n=0,1,2,\dots.
\end{equation*}
It is well known that
\begin{equation}
\sum_{n=0}^\infty C_n\left(\frac14\right)^n=2.
\label{eq:Catalan-gf}
\end{equation}

\begin{lemma} \label{lem:hump-count} 
Let $m\in\NN$.

(i) There are $\binom{2m}{m}$ humps of order $m$.

(ii) There are $C_m$ leading humps of order $m$.
\end{lemma}

\begin{proof}
Each hump of order $m$ corresponds uniquely to a path of $m$ steps starting at $(0,0)$, taking steps $(1,1)$ or $(1,-1)$, and ending at $(2m,0)$. There are $\binom{2m}{m}$ such paths, proving (i). It is well known that exactly $C_m$ of these paths stay on or above the horizontal axis (see Feller \cite[p.~73]{Feller}), which gives (ii). 
\end{proof}

Finally, we will need the following two facts about dyadic rational abscissas.

\begin{lemma} \label{lem:dyadics-are-endpoints}
Let $x$ be a dyadic rational in $(0,1)$. Then

(i) $x$ is an endpoint of some interval $I(x_0)$ with $x_0\in\BB$; and

(ii) the level set $L(T(x))$ is infinite.
\end{lemma}

\begin{proof}
Let $x\in(0,1)$ be a dyadic rational with terminating binary expansion $x=0.\eps_1\eps_2\dots \eps_n$, so that $\eps_n=1$. Let $j$ be the number of $0$'s among $\eps_1,\dots,\eps_n$. If $j<n-j$, we can simply write $x$ as $x=0.\eps_1\dots\eps_n 0^{n-2j}$ and see that $x$ is balanced of order $n-j$, so $x$ is the left endpoint of $I(x)$. In this case, the right endpoint of $I(x)$ is $x':=0.\eps_1\dots\eps_n 0^{n-2j-1}1$. If $j\geq n-j$, put $x'=0.\eps_1\dots\eps_{n-1}01^{2j+2-n}$. Then $x'$ is balanced of order $j+1$, and $x$ is the right endpoint of $I(x')$. This proves (i). Note that in both cases, we have constructed a dyadic rational point $x'$ such that $T(x')=T(x)$ (by Lemma \ref{lem:similar-copies}, since $x$ and $x'$ are the two endpoints of an interval $I(x_0)$ for some $x_0\in\BB$), and the last ``1" in the binary expansion of $x'$ occurs later than the last ``1" in the binary expansion of $x$. Since $x$ was arbitrary, this implies (ii).
\end{proof}

We will see in Section \ref{sec:Baire} that $L(T(x))$ is in fact countable when $x$ is dyadic.

\begin{corollary} \label{cor:all-dyadics}
Every dyadic rational in $(0,1)$ is contained in  $\bigcup_{x_0\in\BB_1}I(x_0)$.
\end{corollary}

\begin{proof}
This follows from Lemma \ref{lem:dyadics-are-endpoints} since each $I(x_0)$ with $x_0\in\BB$ is contained in some $I(x_1)$ with $x_1\in\BB_1$.
\end{proof}

\section{Probability view} \label{sec:probability}

In this section we give more direct proofs of two known results, namely that almost every level set (with respect to Lebesgue measure) is finite, but that the expected cardinality of a level set chosen at random is infinite. We also give a conceptually easier proof of a recent result by Lagarias and Maddock \cite{LagMad1} concerning local level sets. It is natural to think about these results in a probabilistic setting, so we define a probability measure $\rP$ on the range $[0,\frac23]$ by
\begin{equation*}
\rP(A)=\frac32\lambda(A), \qquad A\subset [0,\tfrac23],
\end{equation*}
where $A$ ranges over the Lebesgue subsets of $[0,\frac23]$, and $\lambda$ denotes Lebesgue measure on the line. We let $\rE$ denote the corresponding expectation operator; that is, $\rE(N)=\frac32\int_0^{2/3}N(y)\,dy$ for a Lebesgue measurable function $N:[0,\frac23]\to\RR\cup\{\infty\}$. 

\begin{lemma} \label{lem:leading-humps}
(i) For every hump $H$ there is a leading hump $H'$ of the same order and generation as $H$, such that $\pi_Y(H)=\pi_Y(H')$.

(ii) For every leading hump $H'$, there are only finitely many humps $H$ such that $\pi_Y(H)=\pi_Y(H')$.
\end{lemma}

\begin{proof}
(i) Let $H=H(x_0)$, where $x_0$ is a balanced dyadic rational of order $m$. There is a unique balanced dyadic rational $x_1$ of order $m$  such that $D_j(x_1)=|D_j(x_0)|$ for all $j\leq 2m$. By Lemma \ref{lem:equivalence}, $T(x_0)=T(x_1)$. By definition, $H':=H(x_1)$ is a leading hump of order $m$ and of the same generation as $H$. Hence $H'$ is the same size as $H$ and sits at the same height in the graph of $T$. Therefore, $H'$ is as required.

(ii) This is immediate from Lemma \ref{lem:hump-count}.
\end{proof}

Recall that $\BB_1\subset \BB$ denotes the set of first-generation balanced dyadic rationals. Define a subset $X^*$ of $[0,1]$ by
\begin{equation}
X^*:=[0,1]\backslash\bigcup_{x_0\in\BB_1}I(x_0).
\label{eq:Xstar}
\end{equation}
In other words, $X^*$ is obtained by removing all the dyadic closed intervals above which the graph of $T$ has a first-generation hump. By Corollary \ref{cor:all-dyadics}, $X^*$ does not contain any dyadic rationals other than $0$ and $1$. Note that $X^*$ is symmetric, so the restriction of $T$ to $X^*$ is symmetric as well. The importance of $X^*$ is made clear by the following proposition, which will be used repeatedly throughout this paper.

\begin{proposition} \label{lem:bijection}
The Takagi function $T$ maps $X^*$ onto $[0,\frac12]$. Moreover, $T$ is strictly increasing on $X^*\cap[0,\frac12)$.
\end{proposition}

The proof of the proposition uses the following auxiliary functions. Let
\begin{equation}
T^*(x):=\begin{cases}
T(x), & \mbox{if $x\in X^*$},\\
T(x_0), & \mbox{if $x\in I(x_0)$, where $x_0\in\BB_1$},
\end{cases}
\label{eq:T-star}
\end{equation}
and define piecewise linear approximants of $T^*$ by
\begin{equation}
T_n^*(x):=\begin{cases}
T(x_0), & \mbox{if $x\in I(x_0)$ with $x_0=k/2^{2m}\in\BB_1$ and $2m\leq n$},\\
T_n(x), & \mbox{otherwise}.
\end{cases}
\label{eq:Tn-star}
\end{equation}
Thus, $T_n^*$ permanently ``fixes" each first-generation flat segment in the piecewise linear approximations of $T$.

\begin{lemma} \label{lem:T-star}
The functions $T^*$ and $T_n^*$ have the following properties:

(i) $T_n^*$ is continuous and nondecreasing on $[0,\frac12]$ for every $n$;

(ii) $T_n^*\to T^*$ uniformly on $[0,1]$;

(iii) $T^*$ is continuous and nondecreasing on $[0,\frac12]$.
\end{lemma}

\begin{proof}
Statement (i) is obvious for $n=1$, as $T_1^*=T_1$. In the transition from $T_n^*$ to $T_{n+1}^*$, each horizontal line segment stays fixed, and each line segment of strictly positive integer slope $m$ is replaced with two connecting line segments of slopes $m+1$ and $m-1$, respectively, which meet the original line segment at its endpoints. Thus, it follows inductively that (i) holds for every $n\in\NN$.

Statement (ii) follows from the uniform convergence of $T_n$ to $T$ by the following argument. Let $x\in[0,1]$. If $x\in X^*$, then $T^*(x)-T_n^*(x)=T(x)-T_n(x)$. Suppose instead that $x\in I(x_0)$, with $x_0=k/2^l \in\BB_1$. If $l\leq n$, then the flat segment above $I(x_0)$ has already been fixed by $T_n^*$, and so $T^*(x)-T_n^*(x)=T(x_0)-T(x_0)=0$. Finally, if $l>n$, then 
\begin{equation*}
|T^*(x)-T_n^*(x)|=|T(x_0)-T_n(x)|\leq |T(x_0)-T(x)|+|T(x)-T_n(x)|.
\end{equation*}
Since $|x-x_0|\leq 2^{-l}<2^{-n}$ in this case, the first term on the right hand side of the above inequality can be made uniformly small by choosing $n$ sufficiently large, using the uniform continuity of $T$. Since $T_n\to T$ uniformly, an examination of the three cases above allows us to conclude statement (ii).

Statement (iii), of course, is a direct consequence of (i) and (ii).
\end{proof}

\begin{proof}[Proof of Proposition \ref{lem:bijection}]
We show first that $T$ maps $X^*$ onto $[0,\frac12]$. We have $T^*(0)=0$ and $T^*(\frac12)=T(\frac14)=\frac12$, so Lemma \ref{lem:T-star}(iii) implies that $T^*([0,\frac12])=[0,\frac12]$. Let an ordinate $y\in[0,\frac12]$ be given, and set $x^*:=\min\{x\geq 0:T^*(x)=y\}$. (The minimum is well defined by the continuity of $T^*$.) Suppose, by way of contradiction, that $x^*\in I(x_0)$ for some $x_0\in\BB_1$. Then $T(x_0)=T^*(x^*)=y$. Since $x_0$ is balanced of generation 1 and lies in $[0,\frac12]$, its binary expansion is of the form $x_0=0.\eps_1\eps_2\dots\eps_{2m-1}1$ for some $m$. Let $x_1=0.\eps_1\eps_2\dots\eps_{2m-1}011$. Then $x_1$ is balanced, and the right endpoint of $I(x_1)$ is $x_0$. Thus, $T(x_1)=T(x_0)$. But $T(x_0)=y$ from above, so this contradicts the definition of $x^*$ since $x_1<x_0\leq x^*$. We conclude therefore that $x^*\in X^*$, and as a result, $T(x^*)=T^*(x^*)=y$, as desired.

Next, we show that $T$ is strictly increasing on $X^*\cap[0,\frac12)$. Observe first that, if $x\in X^*\cap(0,\frac12)$, then $D_j(x)>0$ for all $j$. This is true since $x<\frac12$ means that $D_1(x)=1$, and $D_j(x)$ changes by $\pm 1$ at each step, so it could not become negative without first becoming zero, in which case $x$ would not lie in $X^*$. 

Let $x,x'\in X^*\cap[0,\frac12)$ with $x<x'$. Take any dyadic rational $z\in(x,x')$. By Corollary \ref{cor:all-dyadics} there is $x_0\in\BB_1$, say $x_0=k/2^{2m}$, such that $z\in I(x_0)$, and at the same time, neither $x$ nor $x'$ lies in $I(x_0)$. Thus, $I(x_0)\subset (x,x')$. Let $n\geq 2m$. Since $x'$ is not a dyadic rational, the derivative of $T_n^*$ at $x'$ is defined and is equal to the strictly positive number $D_n(x')$. Since $T_n^*$ is nondecreasing, it follows that
\begin{equation}
T(x')\geq T_n(x')=T_n^*(x')>T_n^*(x_0)\geq T_n^*(x)=T_n(x). 
\label{eq:inequality-chain}
\end{equation}
Since $T_n^*(x_0)=T(x_0)$, this implies $T(x')>T(x_0)\geq T_n(x)$.
Letting $n\to\infty$, it follows that $T(x')>T(x)$.
\end{proof}

\begin{remark}
{\rm
(a) The set $X^*\cap[0,\frac12)$ is a proper subset of the set $\frac12\Omega^L$ of Lagarias and Maddock \cite{LagMad2}: $\frac12\Omega^L$ contains certain dyadic rationals (precisely, right endpoints of intervals $I(x_0)$), whereas $X^*\cap[0,\frac12)$ does not. It is exactly this difference which makes $T$ strictly increasing on $X^*\cap[0,\frac12)$. The small price to pay for this is that $X^*$ has higher set-theoretic complexity: whereas $\frac12\Omega^L$ is a closed set, $X^*$ is merely a $G_\delta$.

In fact, a bit more can be said about the structure of $X^*$. It can be gleaned from the first half of the proof of Proposition \ref{lem:bijection} that $X^*\cap[0,\frac12]$ is really the complement of a countable union of half-open intervals $(a,b]$. This is explained by the fact that at each level $T(x_0)$, $x_0\in\BB_1$, countably many of the intervals $I(x)$ glue together at their endpoints to form a half-open interval, and $T$ takes the same value $T(x_0)$ at the endpoints of all these intervals. For instance, the point $x\in X^*$ such that $T(x)=\frac12$ is $x=\frac16$. The removed half-open interval $(\frac16,\frac12]$ is the union of the first-generation intervals $[\frac14,\frac12]$, $[\frac{3}{16},\frac14]$, $[\frac{11}{64},\frac{3}{16}]$, etc., where the left endpoint of the $k$th interval is $\frac12-\sum_{j=1}^k {(\frac14)}^k$.

(b) The function $T^*$ is closely related to the {\em Takagi singular function} $\tau^S$ of Lagarias and Maddock \cite{LagMad1,LagMad2}, as follows: For all $x\in[0,\frac12]$, $T^*(x)=\frac12\tau^S(2x)$. The author is grateful to the referee for pointing this out. Since we do not need this fact here, we omit the details. For more properties of the function $\tau^S$ and the associated singular probability measure, see \cite{LagMad1,LagMad2}.
}
\end{remark}

\begin{definition}
{\rm
Let $\mathcal{G}_T^*=\{(x,T(x)):x\in X^*\}$, so $\mathcal{G}_T^*$ is the graph of $T$ with all first-generation humps removed. For a hump $H=H(x_0)$, the affine transformation that maps $\mathcal{G}_T$ onto $H$ maps $\mathcal{G}_T^*$ onto a set which we call a {\em truncated hump}, and denote by $H^t=H^t(x_0)$. Let $J^t(x_0)=\pi_Y(H^t(x_0))$. Proposition \ref{lem:bijection} implies that if $J(x_0)=[a,a+\frac23{(\frac14)}^m]$, then $J^t(x_0)=[a,a+\frac12{(\frac14)}^m]$. 
}
\end{definition}

In what follows, let $l_y$ denote the horizontal line at level $y$. That is,
\begin{equation*}
l_y:=\{(x,y): x\in\RR\}.
\end{equation*}

\begin{lemma} \label{lem:when-finite}
Let $y\in[0,\frac23]$. Then $|L(y)|<\infty$ if and only if $l_y$ intersects only finitely many leading humps.
\end{lemma}

\begin{proof}
Suppose $l_y$ intersects only finitely many leading humps. Then by Lemma \ref{lem:leading-humps} it intersects only finitely many humps. Let $N$ denote the maximal generation of all humps intersected by $l_y$. Then $l_y$ intersects each hump of generation $N$ only in the corresponding truncated hump, and hence it intersects each such hump in at most two points, by Proposition \ref{lem:bijection}. Thus, $l_y$ has only finitely many intersection points with humps of generation $N$. But then $l_y$ intersects each hump $H$ of generation $N-1$ in at most finitely many points: at most two contributed by the truncated hump $H^t$, plus (at most) finitely many contributed by the humps of generation $N$ contained in $H$. Hence, since $l_y$ intersects only finitely many humps, it has only finitely many intersection points with humps of generation $N-1$. Continuing this argument inductively we find that for each $n\leq N$, $l_y$ has only finitely many intersection points with humps of generation $n$. Since $\mathcal{G}_T$ is a hump of generation $0$, it follows that $|L(y)|<\infty$.

Conversely, suppose $l_y$ intersects infinitely many leading humps, and consider two cases. If there is a number $N$ such that $l_y$ intersects infinitely many leading humps of generation $N$, then $L(y)$ is infinite because leading humps of the same generation can have at most one point in common, and $l_y$ intersects each hump in at least two points. If there is no such $N$, then $l_y$ intersects at least one hump of each generation $n$. If $(x,y)$ is such an intersection point, then there are $n$ indices $k\in\NN$ such that $D_k(x)=0$. It follows that there are $2^{n+1}$ points $x'$ with $|D_k(x)|=|D_k(x')|$ for each $k\in\NN$, so by Lemma \ref{lem:equivalence}, $l_y$ intersects $\GG_T$ at least $2^{n+1}$ times. Since this is the case for every $n\in\NN$, we conclude that $|L(y)|=\infty$.
\end{proof}

(The ``only if" part of the lemma is not used in this section, but will be needed in Section \ref{sec:Baire}.)

\begin{theorem}[Buczolich, 2008 \cite{Buczolich}] \label{thm:finite-ae}
For almost every $y$ (with respect to Lebesgue measure on $[0,\frac23]$), $L(y)$ is a finite set.
\end{theorem}

\begin{proof}
If $y$ is chosen at random and $H$ is a leading hump of order $m$, the probability that the line $l_y$ intersects $H$ is ${(\frac14)}^m$. This gives
\begin{equation*}
\sum_{H\in\HH'}\sP(y\in \pi_Y(H))=\sum_{m=0}^\infty C_m\left(\frac14\right)^m<\infty
\end{equation*}
by \eqref{eq:Catalan-gf}, because there are $C_m$ leading humps of order $m$. Thus, by the Borel-Cantelli lemma, the probability that $l_y$ intersects infinitely many $H$ in $\HH'$ is zero. Therefore, by Lemma \ref{lem:when-finite}, $L(y)$ is finite with probability $1$.
\end{proof}

\begin{lemma} \label{lem:cardinality-of-L}
For any ordinate $y$, if $L(y)$ is finite, then its cardinality is given by
\begin{equation}
|L(y)|=2\cdot|\{H\in\HH: y\in\pi_Y(H^t)\}|.
\label{eq:exact-cardinality}
\end{equation}
\end{lemma}

\begin{proof}
Two truncated humps can intersect each other in at most one point, and if they do, the common point is the lower left or lower right vertex of at least one of the two truncated humps. In that case, the $x$-coordinate of the common point is clearly a dyadic rational, so that the corresponding level set is infinite by Lemma \ref{lem:dyadics-are-endpoints}(ii).


Let $|L(y)|<\infty$, and let $x$ be an abscissa with $T(x)=y$. It must then be the case that $D_j(x)=0$ for only finitely many $j$ (perhaps none), for otherwise, as remarked below Lemma \ref{lem:equivalence}, $L(y)$ would contain a Cantor set. Thus, there is $m\in\ZZ_+$ such that $D_{2m}(x)=0$, but $D_j(x)\neq 0$ for all $j>2m$. This means $(x,y)\in H^t(x_0)$, where $x_0$ is the dyadic rational $k/2^{2m}$ whose first $2m$ digits coincide with those of $x$. Since $x$ was arbitrary, we conclude that each intersection point of the line $l_y$ with $\GG_T$ must lie on some truncated hump, and by the remark at the beginning of the proof, intersections with different truncated humps are disjoint. By Proposition \ref{lem:bijection}, $l_y$ intersects each truncated hump either not at all or in exactly two points. This gives \eqref{eq:exact-cardinality}.
\end{proof}

\begin{theorem}[Lagarias and Maddock, 2010 \cite{LagMad2}] \label{thm:infinite-expectation}
The expected cardinality of a level set $L(y)$ for $y$ chosen at random from $[0,\frac23]$ is infinite. That is,
\begin{equation*}
\rE|L(y)|=\frac32\int_0^{2/3} |L(y)|\,dy=\infty.
\end{equation*}
\end{theorem}

\begin{proof}
By \eqref{eq:exact-cardinality}, $|L(y)|$ is Lebesgue measurable, so its integral is well defined. By Lemma \ref{lem:cardinality-of-L} and Theorem \ref{thm:finite-ae},
\begin{align*}
\rE|L(y)|&=2\sE|\{H\in\HH: y\in\pi_Y(H^t)\}|=2\sum_{H\in\HH}\sP\left(y\in \pi_Y(H^t)\right)\\
&=2\sum_{m=0}^\infty \frac34\binom{2m}{m}\left(\frac14\right)^m=\infty,
\end{align*}
since $\rP(y\in\pi_Y(H^t))=\frac34 \sP(y\in \pi_Y(H))=\frac34{(\frac14)}^m$ for a balanced dyadic $x_0$ of order $m$, and $\binom{2m}{m}\sim 4^m/\sqrt{\pi m}$ as $m\to\infty$.
\end{proof}

\subsection{Local level sets} \label{subsec:local}

Lagarias and Maddock \cite{LagMad1,LagMad2} introduce the concept of a local level set of the Takagi function. They first define an equivalence relation on $[0,1]$ by
\begin{equation}
x\sim x'\quad \stackrel{\defin}{\Longleftrightarrow}\quad |D_j(x)|=|D_j(x')|\ \mbox{for each $j\in\NN$}.
\label{eq:equivalence-relation}
\end{equation}
Note that by Lemma \ref{lem:equivalence}, $x\sim x'$ implies $T(x)=T(x')$. 
The {\em local level set} containing $x$ is then defined by 
\begin{equation*}
L_x^{loc}:=\{x': x'\sim x\}.
\end{equation*}
Lagarias and Maddock point out that each local level set is either finite or a Cantor set.
One of their results deals with the average number of local level sets contained in a level set chosen at random. Let
$N^{loc}(y)$ denote the number of local level sets contained in $L(y)$. 

\begin{remark}
{\rm
Lagarias and Maddock \cite{LagMad1} define $L_x^{loc}$ slightly differently. They treat each dyadic point $x$ as a pair of separate points $x^+$ and $x^-$ according to the two possible binary expansions of $x$, and by \eqref{eq:equivalence-relation}, $x^+$ and $x^-$ represent different local level sets. Thus, compared to our definition, some local level sets contain extra points of the form $x^-$ in the definition of Lagarias and Maddock, where $x^-$ corresponds to the representation of $x$ ending in all ones. But the difference in definitions does not affect the {\em number} of local level sets contained in any level set, which is all we are concerned with in this paper.
}
\end{remark}

\begin{theorem}[Lagarias and Maddock, 2010 \cite{LagMad1}] \label{thm:local-level-sets}
The expected number of local level sets contained in a level set $L(y)$ with $y$ chosen at random from $[0,\frac23]$ is $\frac32$. More precisely,
\begin{equation*}
\rE[N^{loc}(y)]=\frac32\int_0^{2/3}N^{loc}(y)\,dy=\frac32.
\end{equation*}
\end{theorem}

\begin{proof}
Since $L(y)$ is finite for almost every $y$ with respect to Lebesgue measure, we need only consider local level sets which are finite.

\bigskip
\noindent {\bf Claim:} The number of {\em finite} local level sets contained in $L(y)$ is exactly
\begin{equation*}
|\{H\in\HH':y\in\pi_Y(H^t)\}|.
\end{equation*}
To prove the claim, note that there is an obvious one-to-one correspondence between leading humps and truncated leading humps. We show that for $y\in[0,\frac23]$, there is a bijection between the collection of finite local level sets in $L(y)$ and the collection of truncated leading humps which intersect the line $l_y$. 

Note that we can parametrize finite local level sets by their left-most point. If $L_x^{loc}\subset L(y)$ with left-most point $x$, then their must exist $m\in\ZZ_+$ such that $D_{2m}(x)=0$, but $D_j(x)>0$ for all $j>2m$. (Otherwise, $L_x^{loc}$ would be a Cantor set.) Let $x_0$ be the dyadic rational $k/2^{2m}$ whose first $2m$ binary digits coincide with those of $x$. Then $(x,y)\in H^t(x_0)$ as in the proof of Lemma \ref{lem:cardinality-of-L}. Thus, we have a mapping from finite local level sets in $L(y)$ to truncated leading humps which intersect $l_y$. This mapping is bijective. It is onto, because given a truncated leading hump $H^t$ which intersects $l_y$, let $x$ be the left-most point such that $(x,y)\in H^t\cap l_y$; then the local level set $L_x^{loc}$ gets mapped to $H^t$. The mapping is one-to-one, because any truncated hump $H^t$ can have only two points of intersection with $l_y$, and the abscissas of both points belong to the same local level set.

The Claim follows from the above one-to-one correspondence. Using the Claim we obtain, by a calculation similar to that in the proof of Theorem \ref{thm:infinite-expectation},
\begin{equation*}
\rE[N^{loc}(y)]=\sum_{H\in\HH'}\sP\left(y\in \pi_Y(H^t)\right)
=\frac34\sum_{m=0}^\infty C_m\left(\frac14\right)^m=\frac32.
\end{equation*}
Here the second equality follows by Lemma \ref{lem:hump-count}(ii), and the third by \eqref{eq:Catalan-gf}.
\end{proof}

\begin{remark}
{\rm
Lagarias and Maddock \cite{LagMad1} prove Theorem \ref{thm:local-level-sets} above by introducing a ``flattened Takagi function" and using the coarea formula from the theory of functions of bounded variation. The approach used here is more straightforward and appears more natural.
}
\end{remark}

\section{Baire category view} \label{sec:Baire}

We next investigate the level sets from the point of view of Baire category, and find that the results sharply contrast those of the previous section. In this section, descriptive set theory notation is used to indicate the complexity of certain sets. For the definitions and properties of the Borel hierarchy, see Kechris \cite[Sections 11.B and 22]{Kechris}.

Define the set
$$S_\infty:=\{y\in[0,\tfrac23]: |L(y)|=\infty\}.$$

\begin{lemma} \label{lem:dense-union}
For each $n\in\NN$, the set
\begin{equation*}
E_n:=\bigcup_{x_0\in\BB_n}J(x_0)
\end{equation*}
is dense in $[0,\frac23]$.
\end{lemma}

\begin{proof}
Observe that, considered as a stochastic process on $[0,1]$ with Lebesgue measure, $\{D_n\}$ is a symmetric simple random walk, and therefore returns to $0$ infinitely many times with probability one. This implies that for each $n$ the set $U_n:=\bigcup_{x_0\in\BB_n}I(x_0)$ has full measure in $[0,1]$, so in particular, $U_n$ is dense. Since $T$ is continuous, it follows that $T(U_n)$ is dense in $[0,\frac23]$. But $T(U_n)=E_n$.
\end{proof}

\begin{theorem} \label{thm:topology}
(i) For each $m\in\NN$, the set $S_{\leq m}:=\{y: |L(y)|\leq m\}$ is nowhere dense in $[0,\frac23]$.

(ii) The set $S_\infty$ is residual in $[0,\frac23]$, and is ${\bf\Delta}_3^0$ in the Borel hierarchy.
\end{theorem}

\begin{proof}
Recall from the proof of Lemma \ref{lem:when-finite} that, if the line $l_y$ intersects a hump of generation $n$, then $|L(y)|\geq 2^{n+1}$. Choose an integer $n$ such that $2^{n+1}>m$. Then $y\in S_{\leq m}$ implies that $l_y$ does not intersect any hump of generation $n$. But the union of the projections of $n$-th generation humps is dense in $[0,\frac23]$ by Lemma \ref{lem:dense-union}, so each interval in $[0,\frac23]$ contains as a subinterval some $J(x_0)$ with $x_0\in\BB_n$ which does not intersect $S_{\leq m}$. Hence $S_{\leq m}$ is nowhere dense in $[0,\frac23]$, proving (i).

The first statement of (ii) is an immediate consequence of (i),
since $S_\infty=[0,\frac23]\backslash \bigcup_{m=1}^\infty S_{\leq m}$. As for the second statement, a theorem of Borsuk \cite{Borsuk} says that for {\em any} continuous function $f$, the set $S_\infty(f):=\{y:|f^{-1}(y)|=\infty\}$ is the union of a $G_\delta$ and a countable set. Indeed, Lemma \ref{lem:when-finite} gives the explicit form
\begin{equation*}
S_\infty=\bigcap_{\CC\subset\HH',|\CC|<\infty}\bigcup_{H\in\HH'\backslash C}\pi_Y(H).
\end{equation*}
Replacing each $\pi_Y(H)$ with its interior clearly reduces the set on the right by at most countably many points, and turns it into a $G_\delta$ set. Thus $S_\infty$ is the union of a $G_\delta$ and a countable set, and is therefore of type ${\bf\Delta}^0_3$.
\end{proof}

Note that in view of Theorems \ref{thm:finite-ae} and \ref{thm:topology}, the set $\{y:L(y)<\infty\}$ provides yet another example of a set of first category but of full Lebesgue measure in $[0,\frac23]$.

We now investigate the set $S_\infty$ in more detail. Note that it has the natural decomposition
$S_\infty=S_\infty^{co}\cup S_\infty^{uc}$,
where
\begin{gather*}
S_\infty^{co}:=\{y\in[0,\tfrac23]: L(y)\ \mbox{is countably infinite}\},\\
S_\infty^{uc}:=\{y\in[0,\tfrac23]: L(y)\ \mbox{is uncountably infinite}\}.
\end{gather*}

\begin{lemma} \label{lem:uncountables}
The set $S_\infty^{uc}$ has the representation
\begin{equation*}
S_\infty^{uc}=\bigcap_{n=1}^\infty E_n,
\end{equation*}
where $E_n=\bigcup_{x_0\in\BB_n}J(x_0)$ as in Lemma \ref{lem:dense-union}.
\end{lemma}

\begin{proof}
Let $y\in \bigcap_{n=1}^\infty E_n$. Then for each generation $n$ there is a balanced dyadic rational $x_n\in\BB_n$ such that $y\in J(x_n)$. Since each hump of generation $n$ is contained in a hump of generation $n-1$, we may assume that for each $n$, the binary expansion of $x_n$ is an extension of that of $x_{n-1}$. Let $x=\lim_{n\to\infty}x_n$. Then $D_j(x)=0$ for infinitely many $j$, so $L_x^{loc}$ is a Cantor set and is in particular uncountable. It remains to show that $T(x)=y$. But this follows since a hump of generation $n$ is at least of order $n$ (see Definition \ref{def:humps}), so the intervals $J(x_n)$ shrink to the single point $y$. Furthermore,  $T(x_n)\in J(x_n)$. Hence, $T(x)=\lim_{n\to\infty}T(x_n)=y$.

Conversely, suppose $y\not\in \bigcap_{n=1}^\infty E_n$. Since the sets $E_n$ are nested, this means that there is an index $N$ such that $y\not\in E_n$ when $n>N$. In other words, the line $l_y$ intersects humps of only finitely many generations. The rest of the proof is very similar to the first part of the proof of Lemma \ref{lem:when-finite}, but is included here for definiteness. We claim that $l_y$ intersects each hump in at most countably many points. This is certainly true for each hump $H$ of generation $N$: since $l_y$ does not intersect any smaller humps contained in $H$, it intersects $H$ in at most two points by Proposition \ref{lem:bijection}. Thus for any hump $H$ of generation $N-1$, $l_y$ intersects $H$ in at most countably many points contributed by humps of generation $N$ contained in $H$, plus two points contributed by the truncated hump $H^t$. This implies the claim for all humps of generation $N-1$. Continuing this way, we see inductively that the claim is true for all humps. But since the whole graph $\mathcal{G}_T$ is a hump of generation $0$, the claim immediately implies that $L(y)$ is countable, and so $y\not\in S_\infty^{uc}$.
\end{proof}


It is worth noting that, in view of the above lemma, $S_\infty^{uc}$ is the projection onto the $y$-axis of the irregular $1$-set $S_{\rm irr}$ of Buczolich \cite[Theorem 9]{Buczolich}.

We show next that the set $S_\infty^{uc}$ does not contain any dyadic rationals. For each $x_0\in\BB$, write $J(x_0)=[a(x_0),b(x_0)]$, and let $J^\circ(x_0)=(a(x_0),b(x_0))$ denote the interior of $J(x_0)$.

\begin{lemma} \label{lem:endpoints}
If $y$ is a dyadic rational or has a binary expansion of the form $y=0.\eps_1\eps_2\dots\eps_{2n}{(10)}^\infty$ for some $n\in\NN$, then there are at most finitely many $x_0\in\BB$ such that $y\in J^\circ(x_0)$.
\end{lemma}

\begin{proof}
If $x_0\in\BB$, then it has a binary expansion of the form $x_0=0.b_1 b_2\dots b_{2m}$ for some $m\in\NN$. Since $a(x_0)=T(x_0)$, it follows easily from the definition \eqref{eq:Takagi-def} that $a(x_0)$ has a binary expansion of the form $a(x_0)=0.d_1 d_2\dots d_{2m}$. Then $b(x_0)=a(x_0)+\frac23{(\frac14)}^m$ by Lemma \ref{lem:similar-copies}, so that $b(x_0)$ has binary expansion $b(x_0)=0.d_1 d_2\dots d_{2m}{(10)}^\infty$. Thus, if $y\in(a(x_0),b(x_0))$, we must have
\begin{equation}
0.d_1 d_2\dots d_{2m}<y<0.d_1 d_2\dots d_{2m}{(10)}^\infty,
\label{eq:binary-sandwich}
\end{equation}
so that in particular, the first $2m$ binary digits of $y$ must be $d_1,\dots,d_{2m}$.

If $y$ is dyadic, then $y=0.\eps_1\eps_2\dots\eps_{2n}$ for some $n\in\NN$, and \eqref{eq:binary-sandwich} is possible only if $m<n$. Likewise, if $y$ is of the form $y=0.\eps_1\eps_2\dots\eps_{2n}{(10)}^\infty$, \eqref{eq:binary-sandwich} implies $m<n$. Thus, for fixed $n$, there are only finitely many possible choices for $m$, and for each of those, there are only finitely many points $x_0$.
\end{proof}

\begin{proposition} \label{prop:no-dyadics}
The set $S_\infty^{uc}$ does not contain any dyadic rationals. In other words, if $y$ is dyadic, then $L(y)$ is countable. 
\end{proposition}

\begin{proof}
Let $y$ be a dyadic rational, and suppose by way of contradiction that $y\in S_\infty^{uc}$. By Lemmas \ref{lem:uncountables} and \ref{lem:endpoints}, there must be a number $N$ such that
\begin{equation*}
y\in \bigcap_{n=N}^\infty\bigcup_{x_0\in\BB_n}\{a(x_0),b(x_0)\}.
\end{equation*}
Since $b(x_0)$ is never dyadic, this means that for each $n\geq N$ there is a point $x_n\in\BB_n$ such that $y=a(x_n)$. Let $\check{x}_n$ be the unique balanced dyadic rational of generation $n-1$ such that the binary expansion of $x_n$ extends that of $\check{x}_n$. Then $a(x_n)>a(\check{x}_n)$ for each $n\geq N$, because $a(x)=T(x)>0$ for all $x\in\BB_1$. But then $a(\check{x}_n)<y<b(x_n)\leq b(\check{x}_n)$ for each $n\geq N$, contradicting Lemma \ref{lem:endpoints}.
\end{proof}

We observe that level sets at a dyadic level $y$ can be either finite or  infinite. By Lemma \ref{lem:dyadics-are-endpoints}(ii), $L(y)$ is infinite whenever $y=T(x)$ for some dyadic rational $x\in(0,1)$, and this determines a dense set of ordinates $y$. On the other hand, in a separate paper \cite{Allaart2} we show that there are also infinitely many dyadic levels $y$ with $|L(y)|=2$. Examples are $y=1/8$, $y=3/2^7$, and $y=1/2^8$. 

\begin{theorem} \label{thm:S-uncountable}
The set $S_\infty^{uc}$ is ${\bf\Delta}_3^0$ in the Borel hierarchy, and has the decomposition
$S_\infty^{uc}=E\cup M$,
where $E$ is a dense $G_\delta$ given by
\begin{equation*}
E:=\bigcap_{n=1}^\infty \bigcup_{x_0\in\BB_n} (a(x_0),b(x_0)),
\end{equation*}
and $M$ is a countable set disjoint from $E$ which consists exactly of the local maximum ordinate values of $T$. As a result, the set
$\{y\in[0,\frac23]: L(y)\ \mbox{\rm is countable}\}$
is of the first category.
\end{theorem}

\begin{proof}
The set $E$ is obviously a $G_\delta$, and it is dense in $[0,\frac23]$ by Lemma \ref{lem:dense-union} and Baire's theorem. (Since $E_n$ is a union of intervals, removing the endpoints of these intervals does not ruin denseness.) Let $M:=S_\infty^{uc}\backslash E$. By Lemma \ref{lem:uncountables}, $M$ contains only points of the form $a(x_0)$ or $b(x_0)$. If $x_0\in\BB$, then $a(x_0)\not\in S_\infty^{uc}$ by Proposition \ref{prop:no-dyadics}, so $a(x_0)\not\in M$. On the other hand, $b(x_0)$ is the maximum value of $T$ over $I(x_0)$, so $L(b(x_0))$ has at least the same cardinality as $L(\frac23)$. By the result of Kahane \cite{Kahane}, this means $b(x_0)\in S_\infty^{uc}$. Finally, $b(x_0)\not\in E$ by Lemma \ref{lem:endpoints}. Combining these facts, we conclude that $M=\{b(x_0):x_0\in\BB\}$.
That $S_\infty^{uc}$ is ${\bf\Delta}_3^0$ follows since $E$ is ${\bf\Pi}_2^0$ and $M$ is ${\bf\Sigma}_2^0$, so their union is ${\bf\Delta}_3^0$.
\end{proof}

\begin{corollary} \label{cor:countable-level-sets}
The set $S_\infty^{co}=S_\infty\backslash S_\infty^{uc}$ is ${\bf\Delta}_3^0$ in the Borel hierarchy. It is dense, but it is of first category in $[0,\frac23]$ and of Lebesgue measure zero.
\end{corollary}

\begin{proof}
By Theorem \ref{thm:topology}(ii) and Theorem \ref{thm:S-uncountable}, 
$S_\infty^{co}$ is the difference of two ${\bf\Delta}_3^0$ sets, so it is itself ${\bf\Delta}_3^0$. It is dense in view of Proposition \ref{prop:no-dyadics} and the remark following it, but it is of measure zero by Theorem \ref{thm:finite-ae}, and of the first category by Theorem \ref{thm:S-uncountable}.
\end{proof}

\begin{remark}
{\rm
(a) There remain a number of natural questions about the set $S_\infty^{co}$. For instance, could it in fact be an $F_\sigma$, or stronger still, a countable set? Does it contain any points which are {\em not} images of dyadic rationals? By Lemmas \ref{lem:when-finite} and \ref{lem:uncountables}, $y\in S_\infty^{co}$ if and only if the line $l_y$ intersects infinitely many humps, but does not intersect humps of infinitely many generations. It seems difficult to determine how many points $y$ satisfy these two conditions.

(b) Let $S_{<\infty}$ denote the set of ordinates $y$ with $|L(y)|<\infty$. Then the three sets $S_{<\infty}$, $S_\infty^{co}$ and $S_\infty^{uc}$, which together make up the range of $T$, give three different combinations of results when measured both by Lebesgue measure and Baire category:
$S_{<\infty}$ is of full Lebesgue measure but of first category;
$S_\infty^{co}$ is both of measure zero and of first category; and
$S_\infty^{uc}$ is of measure zero but residual (i.e. large in the sense of category).
These results show a remarkable sort of antisymmetry between these three sets.
}
\end{remark}

To end this section, we answer a question posed by Lagarias and Maddock. They show \cite[Theorem 7.2]{LagMad1} that the set
\begin{equation*}
S_\infty^{loc}:=\{y: L(y)\ \mbox{\rm contains infinitely many different local level sets}\}
\end{equation*}
is dense in $[0,\frac23]$, and ask whether it is countable. That this is far from being the case is shown below.

\begin{theorem} \label{thm:infinite-local-level-sets}
(i) The set $S_\infty^{loc}$ is residual (co-meager) in $[0,\frac23]$.

(ii) The set 
\begin{equation*}
S_\infty^{loc,uc}:=\{y: L(y)\ \mbox{\rm contains {\em uncountably many} different local level sets}\}
\end{equation*}
is dense in $[0,\frac23]$, and intersects any subinterval of $[0,\frac23]$ in a continuum. 
\end{theorem}

The author does not know whether $S_\infty^{loc,uc}$ is residual, or whether it has full Hausdorff dimension $1$. Note that the proof of Lemma \ref{lem:uncountables} shows that every uncountable level set contains an uncountable local level set. So in particular, every level set that contains uncountably many local level sets must have at least one of these local level sets being uncountable.

To prove the theorem, we need the following strengthening of Lemma \ref{lem:dense-union}. 

\begin{lemma} \label{lem:stronger-dense-union}
For each $n\in\NN$, the set $T(\BB_n)$ is dense in $[0,\frac12]$.
\end{lemma}

\begin{proof} 
Recall the function $T^*$ defined in \eqref{eq:T-star} and the fact, as stated in Lemma \ref{lem:T-star}, that $T^*$ is continuous. It was also shown in the proof of Proposition \ref{lem:bijection} that $T^*$ maps $[0,1]$ onto $[0,\frac12]$. As in the proof of Lemma \ref{lem:dense-union}, the set $U_n:=\bigcup_{x_0\in\BB_n}I(x_0)$ is dense in $[0,1]$, and so $T^*(U_n)$ is dense in $[0,\frac12]$. But $T^*(U_n)=T(\BB_n)$.
\end{proof}

\begin{corollary} \label{cor:proper-subintervals}
Let $n\in\NN$. For each finite subset $\CC$ of $\BB_n$ and for each interval $I\subset [0,\frac12]$, $I$ properly contains some $J(x_0)$ with $x_0\in\BB_n\backslash \CC$.
\end{corollary}

\begin{proof}
Let $\CC$ be a finite subset of $\BB_n$, and let $I$ be any subinterval of $[0,\frac12]$. Let $\DD$ be the set of all those $x_0\in\BB_n$ for which $\diam(J(x_0))>(1/3)\diam(I)$. Then $\DD$ is finite, so $T(\BB_n\backslash(\CC \cup \DD))$ is dense in $[0,\frac12]$ by Lemma \ref{lem:stronger-dense-union}, since only finitely many points are removed from $T(\BB_n)$. Hence there is an $x_0\in \BB_n\backslash(\CC \cup \DD)$ such that $T(x_0)$ lies in the lower (left) half of $I$.
But $T(x_0)$ is the lower endpoint of $J(x_0)$, so the interval $J(x_0)$, having length at most $(1/3)\diam(I)$, is properly contained in $I$.
\end{proof}

We observe that the interval $[0,\frac12]$ in Corollary \ref{cor:proper-subintervals} can not be replaced by $[0,\frac23]$. For instance, if $n=1$ and $\mathcal{C}=\{\frac14\}$, then the interval $J(\frac14)=[\frac12,\frac23]$ is removed and there are no further first-generation intervals which overlap the interval $[\frac{7}{12},\frac23]$, say.

\begin{lemma} \label{lem:tree-of-humps}
Let $y$ be an ordinate such that the horizontal line $l_y$ at level $y$ intersects, for each $n\in\NN$, (at least) $2^n$ different leading humps $H_{n,1},\dots,H_{n,2^n}$ of generation $n$ in such a way that for each $j=1,\dots,2^{n-1}$, $H_{n,2j-1}$ and $H_{n,2j}$ are part of $H_{n-1,j}$. Then $L(y)$ contains uncountably many local level sets.
\end{lemma}

\begin{proof}
The hypothesis of the lemma says that the line $l_y$ intersects a binary tree of leading humps, with each level of the tree corresponding to a distinct generation, and each hump being a subset of its predecessor in the tree. As the root of the tree we take the hump $H_{0,1}:=\GG_T$, the graph of $T$. To prove the lemma, it clearly suffices to show that each path through this tree starting at the root $H_{0,1}$ will identify a different local level set in $L(y)$.

Let $x_{n,k}\in \BB_n$ be the balanced dyadic rational such that $H_{n,k}=H(x_{n,k})$. Each path through the tree corresponds to a sequence $\{k_n\}$ such that $k_1=1$, and $k_{n+1}$ is either $2k_n-1$ or $2k_n$ for each $n\in\NN$. Since the binary expansion of $x_{n,k_{n+1}}$ extends that of $x_{n,k_n}$, the limit $x:=\lim_{n\to\infty}x_{n,k_n}$ exists, its binary expansion extends that of each $x_{n,k_n}$, and $T(x)=\lim_{n\to\infty}T(x_{n,k_n})=y$ because the sequence $\{H_{n,k_n}\}$ is nested and $\diam(H_{n,k_n})\to 0$. If we have two such paths, coded by sequences $\{k_n\}$ and $\{k_n'\}$, say, then the sequences $\{x_{n,k_n}\}$ and $\{x_{n,k_n'}\}$ differ at some index $n$. But then their respective limits $x$ and $x'$ represent different local level sets, since $D_j(x)\geq 0$ and $D_j(x')\geq 0$ for every $j$, while $D_j(x)\neq D_j(x')$ for some $j$. So $x$ and $x'$ cannot satisfy \eqref{eq:equivalence-relation} in Section \ref{subsec:local}.
\end{proof}

\begin{proof}[Proof of Theorem \ref{thm:infinite-local-level-sets}]
The idea is to first use Corollary \ref{cor:proper-subintervals} to prove the corresponding statements relative to $[0,\frac12]$, and then use the graph's self-similarity properties to extend the results to $[0,\frac23]$. To this end, define intervals
$I_k:=[y_k^*,y_{k+1}^*)$,
where
\begin{equation*}
y_k^*:=\sum_{i=0}^{k-1}\frac{1}{2^{2i+1}}, \qquad k\in\NN.
\end{equation*}
Also put $I_0:=[0,\frac12)$. Note that $[0,\frac23]=\bigcup_{k=0}^\infty I_k\cup\{\frac23\}$, and for each $k\geq 0$, the mapping $\Psi(y):=\frac14 y+\frac12$ maps $I_k$ onto $I_{k+1}$. By Lemma \ref{lem:similar-copies}, we see that if $y\in I_k$, then $L(\Psi(y))$ contains a $\frac14$-scale affine copy of $L(y)$ (two copies in fact, unless $y=0$; see also the graph of $T$). Thus, if $L(y)$ contains infinitely (uncountably) many local level sets, so does $L(\Psi(y))$. It is therefore enough to prove that $S_\infty^{loc}\cap [0,\frac12]$ is residual in $[0,\frac12]$, and that $S_\infty^{loc,uc}$ intersects each subinterval of $[0,\frac12]$ in a continuum.

Let $y\in[0,\frac12]$, and note that if the line $l_y$ intersects infinitely many first-generation leading humps, then $y\in S_\infty^{loc}$, since the intersection points with different first-generation leading humps represent different local level sets. Thus $S_\infty^{loc}$ contains the set
\begin{equation*}
G:=\bigcap_{\CC\subset\BB_1, |\CC|<\infty}\ \bigcup_{x_0\in\BB_1\backslash \mathcal{C}}J^\circ(x_0),
\end{equation*}
where $J^\circ(x_0)$ denotes the interior of $J(x_0)$. The set $G$ is clearly a $G_\delta$, and is dense in $[0,\frac12]$ by Lemma \ref{lem:stronger-dense-union} and Baire's theorem. Hence, $S_\infty^{loc}\cap [0,\frac12]$ is residual in $[0,\frac12]$. As explained above, this proves (i).

To prove (ii), let an interval $V\subset[0,\frac12]$ be given. We will show that there is a continuum of points $y$ in $V$ satisfying the hypothesis of Lemma \ref{lem:tree-of-humps}. This then extends to all of $[0,\frac23]$ in view of the observation at the beginning of the proof.

By Corollary \ref{cor:proper-subintervals}, we can find a first-generation hump $H_{1,1}$ whose projection $J_{1,1}$ onto the $y$-axis lies inside $V$. By Lemma \ref{lem:leading-humps}, we may take $H_{1,1}$ to be a leading hump. Applying Corollary \ref{cor:proper-subintervals} and Lemma \ref{lem:leading-humps} again in the same way, $J_{1,1}$ properly contains the projection $J_{1,2}$ of some first-generation leading hump $H_{1,2}$ distinct from $H_{1,1}$. Since the containment is proper, there is, again by Corollary \ref{cor:proper-subintervals}, at least one other leading hump $H_{1,2}'$ whose projection $J_{1,2}'$ lies fully inside $J_{1,1}\backslash J_{1,2}$. 

Next, we apply Corollary \ref{cor:proper-subintervals} and Lemma \ref{lem:leading-humps} four times to find successively, in the same way as above, second-generation leading humps $H_{2,1},\dots,H_{2,4}$ with projections $J_{2,1},\dots,J_{2,4}$ such that $J_{2,1}\subset J_{1,2}$, $J_{2,k+1}\subset J_{2,k}$ for $k=1,2,3$, and all the containments are proper. Furthermore, these can be chosen so that $H_{2,1}$ and $H_{2,2}$ are part of $H_{1,1}$, while $H_{2,3}$ and $H_{2,4}$ are part of $H_{1,2}$. 
Proceeding in this manner we can build the rest of the binary tree: at stage $n$, we can successively find leading humps $H_{n,1},\dots,H_{n,2^n}$ of generation $n$ with projections $J_{n,1},\dots,J_{n,2^n}$ such that $J_{n,1}\subset J_{n-1,2^{n-1}}$, $J_{n,k+1}\subset J_{n,k}$ for $k=1,\dots,2^n-1$, and all the containments are proper. Furthermore, we can choose these humps in such a way that for each $j=1,\dots,2^{n-1}$, $H_{n,2j-1}$ and $H_{n,2j}$ are part of $H_{n-1,j}$. Note that at each step, there is an alternative leading hump $H_{n,k}'$ whose projection $J_{n,k}'$ also lies inside $J_{n,k-1}$ (or inside $J_{n-1,2^{n-1}}$ if $k=1$) and is disjoint from $J_{n,k}$.

Now that the tree has been constructed, the intersection of all the intervals $J_{n,k}$, where $n=1,2,\dots$ and $k=1,\dots,2^n$, is a single point $y$, and $y\in S_\infty^{loc,uc}$. Since at any step we could have chosen $J_{n,k}'$ instead of $J_{n,k}$, there is in fact a continuum of such binary trees, and each one gives a different ordinate $y$ in $V$ by the assumption that $J_{n,k}$ and $J_{n,k}'$ are disjoint.  Thus, part (ii) of the theorem follows.
\end{proof}

Recall from Theorem \ref{thm:local-level-sets} that the average number of local level sets in a randomly chosen level set is $\frac32$. Theorem \ref{thm:infinite-local-level-sets}(i) stands in marked contrast to that result, and is a further example of the discrepancy between measure theory and Baire category.

\section{A generalization} \label{sec:general}

\begin{figure}
\begin{center}
\epsfig{file=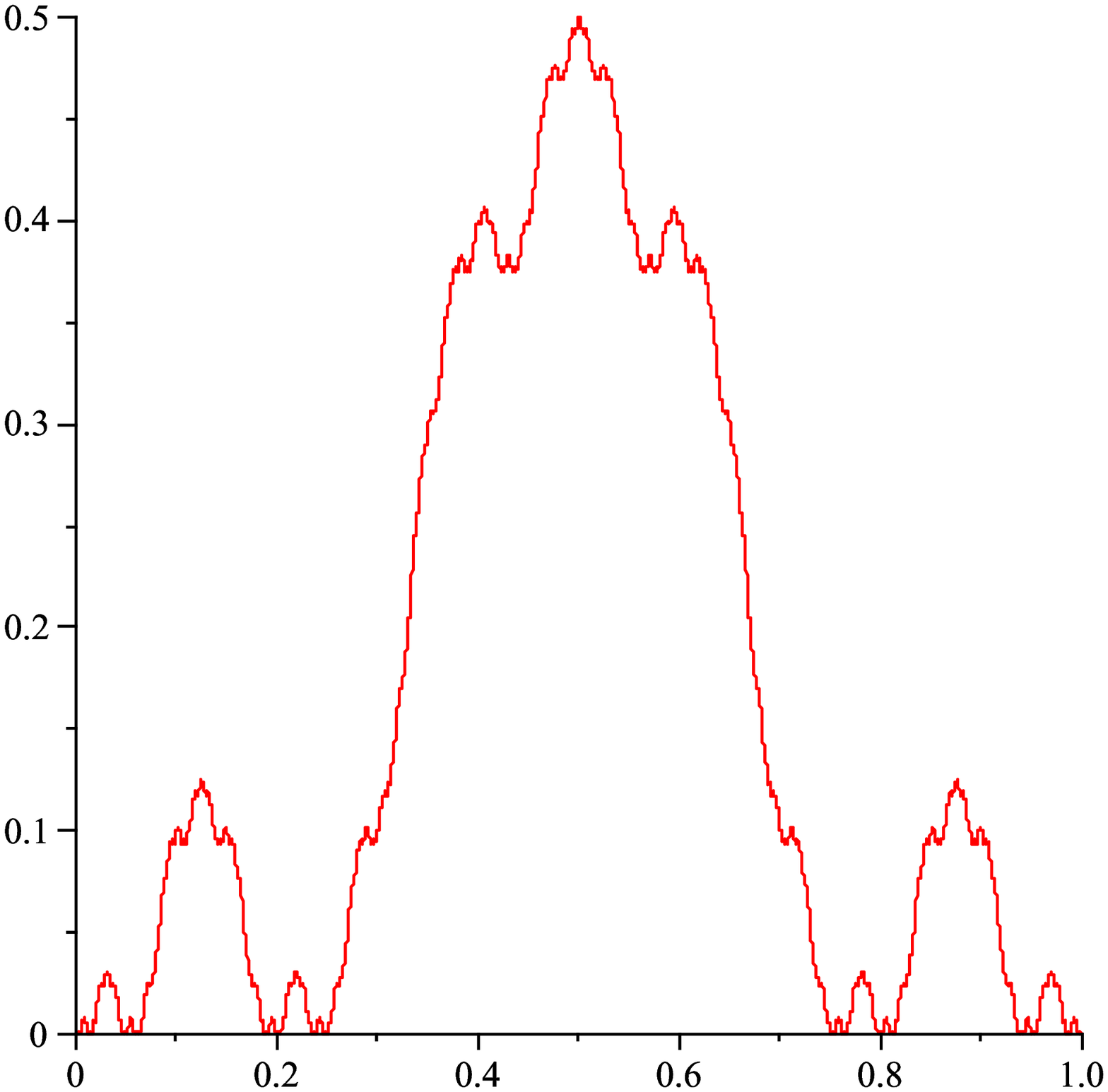, height=.25\textheight, width=.4\textwidth}
\qquad
\epsfig{file=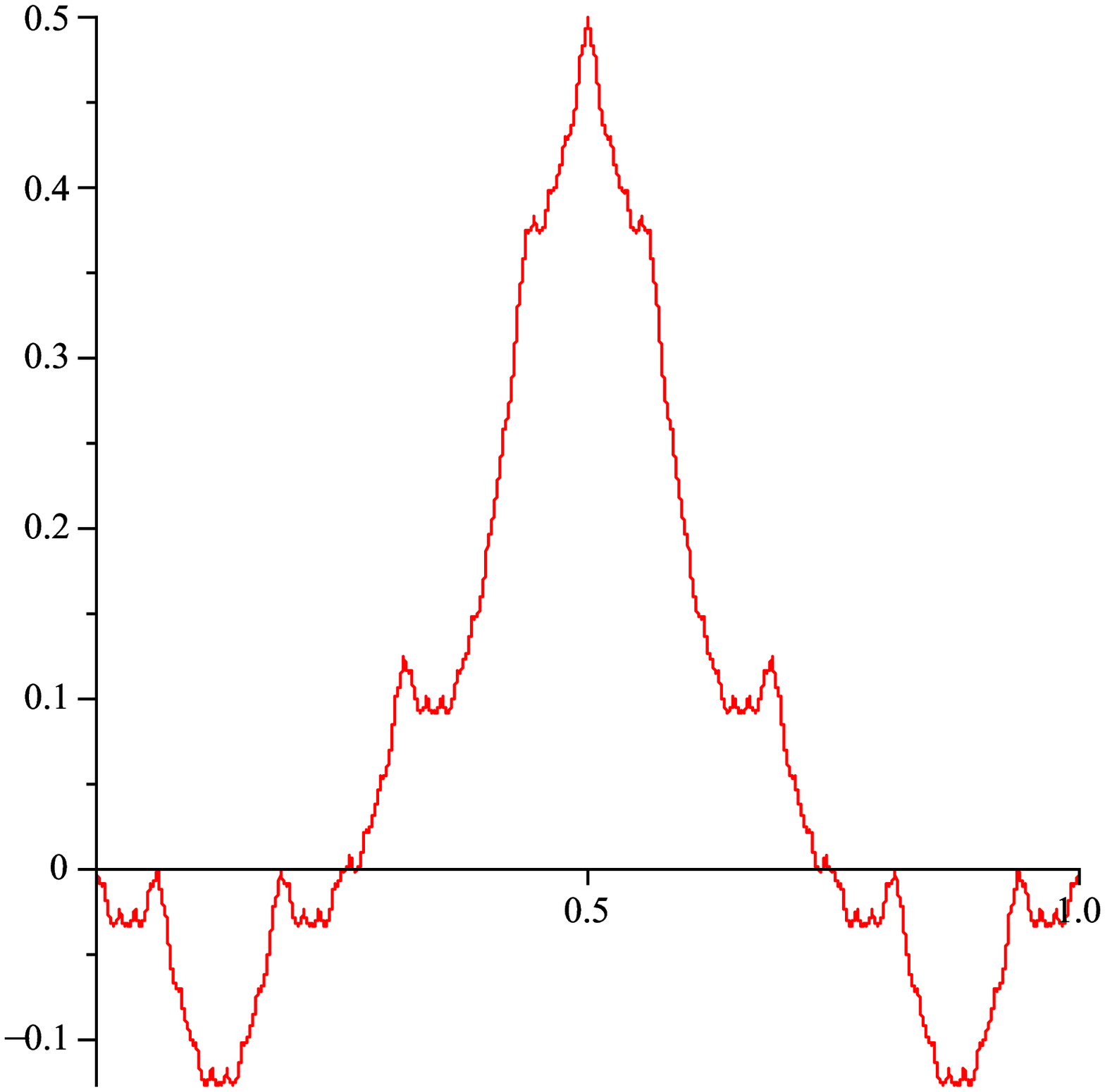, height=.25\textheight, width=.4\textwidth}
\caption{Two functions of the form \eqref{eq:general-function}: the alternating Takagi function (left) has $r_n=(-1)^n$; the function at right has $r_n=1$ if $n\equiv 0\ (\mod 3)$, $r_n=-1$ otherwise.}
\label{fig:general-examples}
\end{center}
\end{figure}

Many of the main results of this paper continue to hold if we multiply the summands in \eqref{eq:Takagi-def} by arbitrary signs. Let ${\bf r}=(r_0,r_1,\dots)$ be a sequence of $\{-1,1\}$-valued numbers, and define
\begin{equation}
f(x)=\sum_{n=0}^\infty \frac{r_n}{2^n}\phi(2^n x).
\label{eq:general-function}
\end{equation}
Two examples are shown in Figure \ref{fig:general-examples}.
Again we let $L(y)=\{x\in[0,1]:f(x)=y\}$ denote the level set of $f$ at level $y$. To analyze these level sets, we need notation analogous to that of Section 2. First, define
\begin{equation*}
f_k(x):=\sum_{n=0}^{k-1}\frac{r_n}{2^n}\phi(2^n x), \qquad k=1,2,\dots.
\end{equation*}
Note that, just as with the Takagi function itself, $f_k\to f$ uniformly on $[0,1]$, and so $f$ is continuous.

With the binary expansion of $x\in[0,1)$ denoted as in \eqref{eq:binary-expansion}, let
\begin{equation*}
D_k(x):=\sum_{j=1}^k r_{j-1}(-1)^{\eps_j}, \qquad k=0,1,\dots.
\end{equation*}
Then $f_k'(x)=D_k(x)$ except at the points $x=j/2^k$, $j=0,\dots,2^k$, where $f_k'(x)$ is undefined. 

\begin{lemma} \label{lem:general-D-expression}
The function $f$ can be expressed in terms of $\{D_n(x)\}$ by
\begin{equation}
f(x)=C({\bf r}) -\frac14\sum_{n=1}^\infty(-1)^{\eps_{n+1}}\frac{D_n(x)}{2^n},
\label{eq:general-D-expression}
\end{equation}
where $C({\bf r})=\sum_{n=0}^\infty r_n/2^{n+2}$ is independent of $x$.
\end{lemma}

\begin{proof}
Let $X_k=(-1)^{\eps_k}$. Then we can write $\phi(x)$ as
\begin{equation*}
\phi(x)=\sum_{k=1}^\infty \frac{\eps_k(1-\eps_1)+(1-\eps_k)\eps_1}{2^k} =\sum_{k=2}^\infty \frac{1-X_1 X_k}{2^{k+1}}.
\end{equation*}
Similarly,
\begin{equation*}
\phi(2^n x)=\sum_{k=2}^\infty \frac{1-X_{n+1} X_{n+k}}{2^{k+1}}.
\end{equation*}
Substituting this into \eqref{eq:general-function} we obtain
\begin{align*}
f(x)&=\sum_{n=0}^\infty \sum_{k=2}^\infty \frac{r_n}{2^{n+k+1}}(1-X_{n+1} X_{n+k})
=\sum_{n=0}^\infty \sum_{j=n+2}^\infty \frac{r_n}{2^{j+1}}(1-X_{n+1}X_j)\\
&=C({\bf r})-\sum_{j=2}^\infty \frac{X_j}{2^{j+1}} \sum_{n=0}^{j-2} r_n X_{n+1}
=C({\bf r})-\sum_{j=2}^\infty \frac{(-1)^{\eps_j}}{2^{j+1}} D_{j-1}(x),
\end{align*}
and re-indexing gives \eqref{eq:general-D-expression}.
\end{proof}

Consequently (and this is of critical importance here), Lemma \ref{lem:equivalence} holds with $f$ replacing $T$. We define balanced dyadic rationals and the notation $\BB$, $\BB_n$ as in Definition \ref{def:balanced}, keeping in mind the new definition of $D_n$. Lemma \ref{lem:similar-copies} no longer holds, of course (except in a few special cases, such as the alternating Takagi function, which arises when ${\bf r}=(1,-1,1,-1,\dots)$). But we will, for consistency, continue to use the word ``hump" for the portion of the graph of $f$ above an interval $[k/2^{2m},(k+1)/2^{2m}]$ if $x_0=k/2^{2m}$ is a balanced dyadic rational. We define the order and generation of a hump as before, and likewise we keep writing $I(x_0)=[k/2^{2m},(k+1)/2^{2m}]$, and $J(x_0)=f(I(x_0))$. Note, however, that $f(x_0)$ need no longer be an endpoint of $J(x_0)$. Different humps may have different shapes, but all are left-to-right symmetric, and we still have ``uniformity within orders": any two humps of the same order $m$ are identical copies of each other. Note also that, while the location of the humps in the graph of $f$ depends on ${\bf r}$, the {\em number} of humps (or leading humps) of a given order does not. Thus -- and this is the second crucial fact -- Lemma \ref{lem:hump-count} continues to hold.

We denote the graph of $f$ by $\GG_f$, and define the {\em height} of the graph by
\begin{equation*}
\height(\GG_f):=\max_{x\in[0,1]}f(x)-\min_{x\in[0,1]}f(x).
\end{equation*}
Similarly, for a hump $H=H(x_0)$ we define its height by
\begin{equation*}
\height(H):=\max_{x\in I(x_0)}f(x)-\min_{x\in I(x_0)}f(x).
\end{equation*}
Since $f(0)=0$ and $f(\frac12)=\pm\frac12$, we have $\height(\GG_f)\geq\frac12$. In case of the Takagi function the height is $\frac23$, and that is also the maximum possible:

\begin{proposition} \label{prop:height}
The height of $\GG_f$ is at most $\frac23$. That is,
\begin{equation*}
\max_{x\in[0,1]}f(x)-\min_{x\in[0,1]}f(x)\leq\frac23.
\end{equation*}
\end{proposition}

\begin{proof}
Set $M:=\max f(x)$ and $m:=\min f(x)$.
Put $s_n:=r_0+r_1+\dots+r_{n-1}$ for $n\in\NN$, and for $j\in\ZZ\backslash\{0\}$, let $\tau_j:=\inf\{n:s_n=j\}$. It is shown in \cite[Theorem 1.1]{Allaart1} that
\begin{equation*}
M=\sum_{k=1}^\infty \left(\frac12\right)^{\tau_{2k-1}},
\end{equation*}
where ${(\frac12)}^\infty$ is interpreted as zero.
(Note that in \cite{Allaart1}, $\phi$ is defined as $2\dist(x,\ZZ)$ and hence the expression given there has an extra factor $2$.) Analogously, we have
\begin{equation*}
m=-\sum_{k=1}^\infty \left(\frac12\right)^{\tau_{1-2k}}.
\end{equation*}
Thus,
\begin{equation*}
M-m=\sum_{j\in\ZZ} \left(\frac12\right)^{\tau_{2j-1}}.
\end{equation*}
Since the $\tau_{2j-1}$ are all distinct and odd, it follows that $M-m\leq\sum_{n=1}^\infty {(\frac12)}^{2n-1}=\frac23$.
\end{proof}

As a consequence, we have for each hump $H$ of order $m$, 
\begin{equation*}
\frac12\left(\frac14\right)^m\leq\height(H)\leq \frac23\left(\frac14\right)^m.
\end{equation*}
This follows since $H$ is itself the (scaled and shifted) graph of a function of the form \eqref{eq:general-function}, with ${\bf r}'=(r_{2m},r_{2m+1},\dots)$.

Next, define a set
\begin{equation*}
X:=[0,1]\backslash \bigcup_{x_0\in\BB_1} I^\circ(x_0),
\end{equation*}
where $I^\circ(x_0)$ denotes the interior of $I(x_0)$. Note that $X$ is closed. (Here we slightly deviate from the way we defined $X^*$ earlier, for reasons that will be made clear below.) For a hump $H(x_0)$ of generation $n$, define the subset $X(x_0)$ of $I(x_0)$ similarly by
\begin{equation*}
X(x_0):=I(x_0)\backslash \bigcup_{x_1\in\BB_{n+1}} I^\circ(x_1).
\end{equation*}
We call the graph of $f$ restricted to $X(x_0)$ a {\em truncated hump}, and denote it by $H^t=H^t(x_0)$. The reason for including the endpoints of the intervals $I(x_1)$ in $X(x_0)$ is to ensure that $\pi_Y(H^t)$ is again an interval.

Since Lemma \ref{lem:equivalence} still holds, Lemma \ref{lem:leading-humps} remains valid as well. In place of the key Proposition \ref{lem:bijection} we have the following result, slightly weaker but sufficient for our purposes.

\begin{lemma} \label{lem:key-lemma-general}
(i) If $r_0=1$, then $f(X)=[0,\frac12]$. If $r_0=-1$, then $f(X)=[-\frac12,0]$.

(ii) If $y\not\in f(\BB_1)$, then $L(y)$ intersects $X\cap[0,\frac12]$ in at most one point.

(iii) For every $y$, $L(y)\cap X$ is countable.
\end{lemma}

\begin{proof}
Assume throughout that $r_0=1$; the case $r_0=-1$ is entirely similar. Define a function
\begin{equation*}
f^*(x):=\begin{cases}
f(x), & \mbox{if $x\in X$},\\
f(x_0), & \mbox{if $x\in I(x_0)$, where $x_0\in\BB_1$}.
\end{cases}
\end{equation*}
Note that $f^*$ is well defined, because if $x$ is an endpoint of an interval $I(x_0)$, we have $f(x)=f(x_0)$. For this reason the use of $X$ instead of $X^*$ as defined in \eqref{eq:Xstar} is immaterial, and in particular, $f^*$ is a generalization of the function $T^*$ from \eqref{eq:T-star}. Define piecewise linear approximants of $f^*$ by
\begin{equation*}
f_n^*(x):=\begin{cases}
f(x_0), & \mbox{if $x\in I(x_0)$ with $x_0=k/2^{2m}\in\BB_1$ and $2m\leq n$},\\
f_n(x), & \mbox{otherwise}.
\end{cases}
\end{equation*}
By the same argument as in the proof of Lemma \ref{lem:T-star}, each $f_n^*$ is continuous and nondecreasing on $[0,\frac12]$, and $f_n^*\to f^*$ uniformly in $[0,1]$. Thus, $f^*$ is continuous and nondecreasing on $[0,\frac12]$ as well. It is easy to see that $f(X)=f^*([0,1])$.
Since $f^*(0)=0$ and $f^*(\frac12)=\frac12$, the symmetry of the graph of $f$ gives $f^*([0,1])=[0,\frac12]$. This proves (i).

The proof of (ii) is very similar to the second part of the proof of Proposition \ref{lem:bijection}. Let $x,x'\in X\cap[0,\frac12]$ with $x<x'$, and assume neither $x$ nor $x'$ is an endpoint of an interval $I(x_0)$ with $x_0\in\BB_1$. By analogy with \eqref{eq:inequality-chain}, we have
\begin{equation*}
f_n(x')=f_n^*(x')\geq f_n^*(x_0)=f(x_0)\geq f_n^*(x)=f_n(x), 
\end{equation*}
for a certain $x_0\in\BB_1$ and $n$ sufficiently large. Letting $n\to\infty$ we obtain $f(x')\geq f(x_0)\geq f(x)$. Thus, if $f(x)=f(x')=y$, it must be the case that $y\in f(\BB_1)$. This proves (ii).

For (iii), we need only consider the case $y=f(x_0)$ with $x_0\in\BB_1$. Points $x\in X$ with $f(x)=y$ can be of two types: (a) endpoints of intervals $I(x_0)$ with $x_0\in\BB_1$, or (b) proper limits of sequences of such endpoints. Clearly, there are only countably many points of type (a). We claim there are at most two points in $[0,\frac12]$ of type (b). To see this, let $x\in X\cap[0,\frac12]$, and suppose there are $x_0,x_1\in\BB_1$ with $x_0<x<x_1$ and $f(x_0)=f(x_1)=y$. Then for all large enough $n$, $f_n^*=y$ everywhere on $I(x_0)$ and $I(x_1)$, and since $f_n^*$ is nondecreasing on $[0,\frac12]$, it follows that $f_n^*$ is constant on $(x_0,x_1)$, an interval containing $x$. But this would mean $x\in I(x_2)$ for some $x_2\in\BB_1$, and being a member of $X$, $x$ must be an endpoint of $I(x_2)$. Thus, the only points of type (b) with $f(x)=y$ lie either to the left of every interval $I(x_0)$ with $f(x_0)=y$, or to the right of each such interval. Hence there can be at most two such points. This establishes (iii).
\end{proof}

We are now ready to restate some of the main results of this paper for the more general function $f$. First, let $[c,d]$ denote the range of $f$, and normalize Lebesgue measure on $[c,d]$ to a probability measure $\rP$ by
$\rP(A)=(d-c)^{-1}\lambda(A)$, $A\subset [c,d]$.
Let $\rE$ denote the corresponding expectation operator. We define local level sets exactly as in Section \ref{subsec:local}, keeping in mind that the configuration of a local level set $L_x^{loc}$ now depends on ${\bf r}$ via the sequence $\{D_n(x)\}$. Let $N^{loc}(y)$ denote the number of local level sets contained in $L(y)$.

\begin{theorem} \label{thm:main-general}
(i) For almost every $y\in[c,d]$, $L(y)$ is finite.

(ii) The expected cardinality of $L(y)$ with $y$ chosen at random from $[c,d]$ is infinite:
\begin{equation*}
\rE|L(y)|=(d-c)^{-1}\int_c^d |L(y)|=\infty.
\end{equation*}

(iii) The average number of local level sets contained in a level set is at least $\frac32$, but at most $2$:
\begin{equation*}
\rE[N^{loc}(y)]=(d-c)^{-1}\int_c^d N^{loc}(y)\,dy \in\left[\tfrac32,2\right].
\end{equation*}

(iv) The set $S_\infty^{uc}(f):=\{y:L(y)\ \mbox{is uncountable}\}$ can be represented as
\begin{equation*}
S_\infty^{uc}(f)=\bigcap_{n=1}^\infty \bigcup_{x_0\in\BB_n}J(x_0),
\end{equation*}
and $S_\infty^{uc}(f)$ is residual in $[c,d]$.
\end{theorem}

The proofs are essentially the same as before, substituting Lemma \ref{lem:key-lemma-general} for Proposition \ref{lem:bijection}. Other lemmas must be slightly modified as well: for instance, Lemma \ref{lem:when-finite} should be replaced by

\begin{lemma} \label{lem:when-finite-general}
Suppose $y\not\in f(\BB)$. Then $|L(y)|<\infty$ if and only if $l_y$ intersects only finitely many leading humps.
\end{lemma}

The same condition on $y$ should be added to Lemma \ref{lem:cardinality-of-L} and to the Claim in the proof of Theorem \ref{thm:local-level-sets}. Since the set $f(\BB)$ is countable it has measure zero, and therefore its exception has no effect on the probabilistic results (i)-(iii) in Theorem \ref{thm:main-general}. The bounds in (iii) follow since for a hump $H$ of order $m$, the height of $H^t$ is exactly $\frac12{(\frac14)}^m$ by Lemma \ref{lem:key-lemma-general}(i), so by Proposition \ref{prop:height}, $\frac34{(\frac14)}^m\leq \rP(y\in\pi_Y(H^t)) \leq {(\frac14)}^m$.
Statement (iv) follows since Lemmas \ref{lem:dense-union} and \ref{lem:uncountables} continue to hold verbatim in the general case, and together they imply that $S_\infty^{uc}(f)$ contains a dense $G_\delta$, so it is residual. (The proof of Lemma \ref{lem:uncountables} requires only a minor modification, namely that $l_y$ intersects each truncated hump in at most {\em countably} many points rather than at most two points. This explains the need for Lemma \ref{lem:key-lemma-general}(iii).)

\begin{remark}
{\rm
We should note that not all of the statements in Sections \ref{sec:probability} and \ref{sec:Baire} transfer to the general case. For instance, Proposition \ref{prop:no-dyadics} is not true in general: If $f$ is the alternating Takagi function
\begin{equation*}
f(x)=\sum_{n=0}^\infty \frac{(-1)^n}{2^n}\phi(2^n x),
\end{equation*}
depicted at left in Figure \ref{fig:general-examples},
then it is not too hard to see that $f\geq 0$, and $L(0)$ is the ``middle half" Cantor set and is hence uncountable. Likewise, $L(f(x_0))$ is uncountable for each $x_0\in\BB$. On the other hand, $L(\frac12)=\{\frac12\}$, so here we have the opposite situation to the case of the Takagi function: $S_\infty^{uc}(f)$ contains the minimum but not the maximum ordinate value of each hump.
}
\end{remark}

\section*{Acknowledgment}
I wish to thank Jeffrey Lagarias and Zachary Maddock for sending preprints of their papers, which were the principal source of inspiration for this article; and David Simmons for pointing out a gap in the original proof of Theorem \ref{thm:finite-ae}. I am greatly indebted to the referee for his or her thorough reading of the paper and for the long list of valuable comments and corrections, especially concerning the proof of Theorem \ref{thm:infinite-local-level-sets}.

\end{document}